# A functional recurrence to obtain the prime numbers using the Smarandache prime function.

## Sebastián Martín Ruiz

Theorem: We are considering the function:

For $n$ integer:

$$F(n) = n + 1 + \sum_{m=n+1}^{2n} \prod_{i=n+1}^{m} \left[ - \left[ - \frac{\sum_{j=1}^{i}\left(\left\lfloor \frac{i}{j} \right\rfloor - \left\lfloor \frac{i-1}{j} \right\rfloor\right) - 2}{i} \right] \right]$$

one has: $p_{k+1} = F(p_k)$ for all $k \geq 1$ where $\{p_k\}_{k \geq 1}$ are the prime numbers and $\lfloor x \rfloor$ is the greatest integer less than or equal to x.

Observe that the knowledge of $p_{k+1}$ only depends on knowledge of $p_k$ and the knowledge of the fore primes is unnecessary.

Proof:

Suppose that we have found a function $P(i)$ with the following property:

$$P(i) = \begin{cases} 1 \text{ if } i \text{ is composite} \\ 0 \text{ if } i \text{ is prime} \end{cases}$$

This function is called Smarandache prime function.(Ref.)

Consider the following product:

$$\prod_{i=p_k+1}^{m} P(i)$$

If $p_k < m < p_{k+1}$ $\prod_{i=p_k+1}^{m} P(i) = 1$ since $i : p_{k+1} \leq i \leq m$ are all composites.

If $m \geq p_{k+1}$ $\prod_{i=p_k+1}^{m} P(i) = 0$ since $P(p_{k+1}) = 0$

Here is the sum:

$$\sum_{m=p_k+1}^{2p_k} \prod_{i=p_k+1}^{m} P(i) = \sum_{m=p_k+1}^{p_{k+1}-1} \prod_{i=p_k+1}^{m} P(i) + \sum_{m=p_{k+1}}^{2p_k} \prod_{i=p_k+1}^{m} P(i) = \sum_{m=p_k+1}^{p_{k+1}-1} 1 =$$

$$= p_{k+1} - 1 - (p_k + 1) + 1 = p_{k+1} - p_k - 1$$

The second sum is zero since all products have the factor $P(p_{k+1}) = 0$.

Therefore we have the following recurrence relation:

$$p_{k+1} = p_k + 1 + \sum_{m=p_k+1}^{2p_k} \prod_{i=p_k+1}^{m} P(i)$$

Let's now see we can find $P(i)$ with the asked property.

Consider:

$$\left\lfloor \frac{i}{j} \right\rfloor - \left\lfloor \frac{i-1}{j} \right\rfloor = \begin{cases} 1 & si \quad j \mid i \\ 0 & si \quad j \text{ not} \mid i \end{cases} \quad j = 1, 2, \cdots, i \quad i \geq 1$$

We deduce of this relation:

$$d(i) = \sum_{j=1}^{i} \left\lfloor \frac{i}{j} \right\rfloor - \left\lfloor \frac{i-1}{j} \right\rfloor$$

where $d(i)$ is the number of divisors of $i$.

If $i$ is prime $d(i) = 2$ therefore:

$$-\left\lfloor -\frac{d(i)-2}{i} \right\rfloor = 0$$

If $i$ is composite $d(i) > 2$ therefore:

$$0 < \frac{d(i)-2}{i} < 1 \Rightarrow -\left\lfloor -\frac{d(i)-2}{i} \right\rfloor = 1$$

Therefore we have obtained the Smarandache Prime Function $P(i)$ which is:

$$P(i) = -\left\lfloor -\frac{\sum_{j=1}^{i}\left(\left\lfloor \frac{i}{j} \right\rfloor - \left\lfloor \frac{i-1}{j} \right\rfloor\right) - 2}{i} \right\rfloor \quad i \geq 2 \quad \text{integer}$$

With this, the theorem is already proved.